\documentclass[11pt]{article}
\usepackage{sectsty,latexsym,amssymb,dsfont,textcomp}
\usepackage[centertags]{amsmath}
\usepackage{fullpage}
\usepackage{color}
\usepackage{epsfig,psfrag}

\usepackage[centertags]{amsmath}
\sectionfont{\normalsize}
\subsectionfont{\normalsize}
\begin{document}
\date{}
\numberwithin{equation}{section}
\title{Relative entropy in multi-phase models of 1d elastodynamics:
     Convergence of a non-local to a local model
}

\author{
Jan Giesselmann\footnote{University of Stuttgart, Stuttgart, Germany,
     jan.giesselmann@mathematik.uni-stuttgart.de}
}

\maketitle


\newtheorem{theorem}{Theorem}[section]
\newtheorem{prop}[theorem]{Proposition}
\newtheorem{lemma}[theorem]{Lemma}
\newtheorem{cor}[theorem]{Corollary}

\newtheorem{definition}[theorem]{Definition}
\newtheorem{notation}[theorem]{Notation}
\newtheorem{example}[theorem]{Example}
\newtheorem{conj}[theorem]{Conjecture}
\newtheorem{prob}[theorem]{Problem}

\newtheorem{remark}[theorem]{Remark}

\newcommand{\bB}[1]{\subsubsection*{{\normalsize Proof #1:} $\vspace*{-0.3cm}$}}
\newcommand{\eB}{\hfill{$\blacksquare$}\bigskip}

\newcommand{\wkarr}{\; \rightharpoonup \;}
\def\Weak{\,\,\relbar\joinrel\rightharpoonup\,\,}

\font\msym=msbm10
\def\charf {\mbox{{\text 1}\kern-.24em {\text l}}} 

\def\A{\mathbb A}
\def\Z{\mathbb Z}
\def\K{\mathbb K}
\def\J{\mathbb J}
\def\L{\mathbb L}
\def\D{\mathbb D}
\def\cD{\mathcal D}
\def\cO{\mathcal O}
\def\cQ{\mathcal Q}
\def\Mink{{\mathop{\hbox{\msym \char '115}}}}
\def\Integers{{\mathop{\hbox{\msym \char '132}}}}
\def\Complex{{\mathop{\hbox{\msym\char'103}}}}
\def\C{\Complex}
\font\smallmsym=msbm7

\newcommand{\del}{\partial}
\newcommand{\eps}{\varepsilon}
\newcommand{\cof}{\hbox{cof}\,}
\newcommand{\dt}{\operatorname{d}t}
\newcommand{\bx}{{\bf x}}
\newcommand{\bn}{{\bf n}}
\newcommand{\bv}{{\bf v}}
\newcommand{\tv}{{\tilde v^\alpha}}
\newcommand{\btv}{{\bf {\tilde v}}^\alpha}
\newcommand{\bF}{{\bf F}}
\newcommand{\btF}{{\bf \tilde F}^\alpha}
\newcommand{\btC}{{\bf \tilde C}^\alpha}
\newcommand{\tvi}{\tilde v_i^\alpha}
\newcommand{\tfij}{\tilde F_{ij}^\alpha}
\newcommand{\tcij}{\tilde C_{ij}^\alpha}
\newcommand{\fij}{F_{ij}}
\newcommand{\dwdf}{\frac{\del W}{\del \fij}}
\newcommand{\dwdff}{\frac{\del^2 W}{\del \fij \del F_{kl}}}
\newcommand{\dx}{\operatorname{d}x}
\newcommand{\dy}{\operatorname{d}y}
\newcommand{\dxx}{\operatorname{d}\bx}
\newcommand{\ds}{\operatorname{d}{\bf s}}
\newcommand{\ddt}{\frac{\operatorname{d}}{\operatorname{d}t}}

\newcommand{\tic}{\tilde c^\alpha}
\newcommand{\tiw}{\tilde w^\alpha}
\newcommand{\tiu}{\tilde u^\alpha}

\newcommand{\Id}{{\operatorname{Id}}}
\newcommand{\cP}{\mathcal{P}}
\newcommand{\cS}{\mathcal{S}}
\newcommand{\cM}{\mathcal{M}}
\renewcommand{\div}{\hbox{div}\,}

\newcommand{\setN}{\mathbb{N}}
\newcommand{\setR}{\mathbb{R}}
\def\Real{{\mathbb{R}}}
\def\R{\Real}
\def\torus{{\mathbb{T}}}
\def\T{\torus}
\def\charf {{{\text{\rm 1}}\kern-.24em {\text{\rm l}}}}

\newcommand{\BBR}[1]{\left( #1 \right)}

\def\div{\hbox{div}\,}
\def\supp{\hbox{supp}\,}
\def\dist{\hbox{dist}\,}
\newcommand{\tcb}{\textcolor{blue}}
\newcommand{\tcg}{\textcolor{green}}
\newcommand{\tcr}{\textcolor{red}}

%
%

%

\begin{abstract}
\noindent
In this paper we study a local and a non-local regularization of the system of nonlinear elastodynamics with a non-convex energy.
We show that solutions of the non-local model converge to those of the local model in a certain regime.
The arguments are based on the relative entropy framework and provide an example how local and non-local regularizations may compensate 
for non-convexity of the energy and enable the use of the relative entropy stability theory -- even if the energy is not quasi- or
poly-convex.
\end{abstract}

%
%
%
%

\section{Introduction}
This paper is concerned with the relation between different models for  shearing motions of an elastic bar in one space dimension.
The models under consideration are based on the equations of nonlinear elastodynamics.
As a multi-phase situation is to be described the energy density is a non-convex function, such that the first order problem without viscosity and capillarity is of  hyperbolic-elliptic type, \cite{Jam80,Tru93}.
In this situation entropy conditions - which are standard in the study of hyperbolic conservation laws - are not sufficient to guarantee uniqueness of weak solutions, see \cite{AK91,LeF02}.
The same difficulty is present in the description of phase transitions in compressible fluids.

A classical way to resolve this problem is the introduction of higher order regularizing terms, which, in addition, can be understood as modeling surface tension.
This strategy goes back to the works of van der Waals and Korteweg in the late 19th century \cite[e.g]{Kor01,vdW88}.
In the last decades {\it local}, i.e., second (deformation) gradient, regularizations were considered by many authors, see \cite[e.g.]{AB82,CL01,JTB02,Sle83}.

However, {\it non-local} regularizations, see \cite{Roh05}, involving convolution terms instead of higher order derivatives, have some advantages from a statistical 
mechanics viewpoint \cite{ABCP96,FM98,RT97}.
Therefore, it is important to understand the relations between both classes of models.

In most studies on non-local models arguments based on Taylor expansions were presented indicating that the regularization term of the non-local model converges to the one 
in the local model.
However, these arguments are purely formal, assuming an amount of regularity not guaranteed by the equations, and only consider convergence of one of the operators in the equation.

We will show that solutions of a particular class of non-local models \eqref{eq:non} indeed converge to solutions of a local model \eqref{eq:loc} in  a scaling which avoids the sharp 
interface limit.
In doing so we provide an example how local and non-local regularizations enable the use of a modified relative entropy framework for the derivation of stability results,
 in case of entropies
which are not poly- or quasi-convex. This complements the  results found in \cite{Gie14}.
In addition, we use the modified relative entropy framework to provide an easy argument showing that solutions of the non-local model continuously depend on initial data.

The relative entropy framework was introduced for hyperbolic problems with convex entropy in \cite{Daf79,Dip79}.
In recent years it was successfully used in the study of hyperbolic conservation laws and related systems.
For a general overview on its development in the last decades and its extension to quasi- and poly-convex entropies we refer the reader to the references in \cite[Section 5.7]{Daf10}.
Recent results based on the relative entropy technique include \cite[e.g.]{BT13,FJN12,FN12,JJN13,LT13,LV11}.

Let us describe the models under consideration in more detail.
We consider the following local one-dimensional nonlinear elasticity  model:
\begin{equation}\label{eq:loc}
 \begin{split}
  u_t - v_x &=0\\
  v_t -W'(u)_x &= \mu v_{xx} -\gamma u_{xxx},
 \end{split}
\end{equation}
where $W \in C^3(\setR,[0,\infty))$ is the (non-convex) energy density, $\mu>0$ is a viscosity parameter and $\gamma>0$ is a capillarity parameter.
Our only assumption on $W,$ apart from its regularity and non-negativity, is that there exists some $\bar W>0$ such that
\begin{equation}\label{ass:W}
 W''(u) > - \bar W \quad \forall u \in \setR.
\end{equation}
This allows for very general multi-well structures of $W.$
The fact that $W$ is defined on all of $\setR$ is due to our focus on shearing motions. If we considered \eqref{eq:loc} as a model for longitudinal motions 
$W$ would only be defined for $u>0$. This is elaborated upon in Remark \ref{rem:state}.
The model \eqref{eq:loc} -- with $W$ being a double-well potential -- was considered as a model for compressible liquid-vapor flows
 as well as shearing motions in an elastic bar in  \cite[e.g]{AB82,CL01,HL00,Sle83,Sle84}.
 
We are going to compare \eqref{eq:loc} to a family of non-local models parametrized by $\eps>0$
\begin{equation}\label{eq:non}
 \begin{split}
  u^\eps_t - v^\eps_x &=0\\
  v^\eps_t -W'(u^\eps)_x &= \mu v^\eps_{xx} -L_\eps[u^\eps]_x,
 \end{split}
\end{equation}
where 
\begin{equation}\label{nlo} L_\eps[u]:=\frac{1}{\eps^2} \big( \phi_\eps * u -u\big),\end{equation}
$*$ denotes convolution, and $\phi$ is some mollifier satisfying
\begin{equation}\label{prop:phi}\begin{split}
& \phi \in C^\infty_{0}(\setR,[0,\infty)),\ \operatorname{supp}(\phi) \subset [-1,1],\ \phi(x)=\phi(-x) \ \forall x \in \setR,\\
& \int \phi(x)\dx =1, \ \int \phi(x) x^2\dx =2\gamma,
\end{split}\end{equation}
and $\phi_\eps(\cdot):=\tfrac{1}{\eps}\phi(\cdot/\eps).$
For the modeling background of \eqref{eq:non}, \eqref{nlo}, see \cite{Roh05_IFB}. Well-posedness analysis for \eqref{eq:non} considering weak solutions can be found in \cite{DR08}
 while numerical results were presented in
 \cite{HR08,Kre10} and sharp interface limits were investigated in \cite{Roh05_IFB}.
In addition, the convergence of local to non-local models was investigated numerically in \cite{Kre10}.
 In the experiments presented there convergence is observed but no convergence rates were determined.

The lower order model in \cite{Gie14} can be seen as a special case of \eqref{eq:non}, \eqref{nlo} with the convolution kernel
being the Green's function of a screened Poisson equation.

The related problem of convergence of non-local to local Navier-Stokes-Korteweg models was investigated in \cite{CH11,CH13} for a particular convolution kernel
 and densities close to a constant state $\bar \rho$ 
satisfying $W''(\bar \rho)>0.$ 
The approach in \cite{CH11,CH13}   is based on the analysis of Fourier modes.
 An extension of that approach removing the assumptions
 on the initial density can be found in \cite{Cha14}.

The factor $\eps^{-2}$ in front of the non-local term in \eqref{eq:non}, \eqref{nlo} might be surprising. However, it can also be found in \cite{Cha14,CH11,CH13} and
   is indeed necessary in order to decouple the 
non-local to local limit from the sharp interface limit.

In order to avoid problems introduced by boundary conditions we consider  \eqref{eq:loc} and \eqref{eq:non} on the flat  unit circle $S^1$, that is to say 
the unit interval with periodic boundary conditions.
In particular, this precludes problems when defining the convolution operator $L_\eps$ near the boundary.

The remainder of this paper is organized as follows. We recall energy balances related to \eqref{eq:loc} and \eqref{eq:non} in Section   \ref{sec:energy}.
Section \ref{sec:wp} is devoted to establishing the well-posedness of \eqref{eq:loc} and \eqref{eq:non}.
The non-local to local limit is derived using a relative entropy estimate in Section \ref{sec:limit}.

\subsection{Energy estimates}\label{sec:energy}

Let us recall the well-known fact that strong solutions of \eqref{eq:loc} and \eqref{eq:non} satisfy energy balance laws
 which are compatible with the principle of material frame indifference.
To keep this paper self-contained we prove the two subsequent lemmas which can be found in many places in the literature.
\begin{lemma}[Energy balance of the local model]
 Let $T,\mu,\gamma>0$ be given. Let $(u,v)$ be a strong solution of \eqref{eq:loc} in $(0,T)\times S^1.$
Then, in $(0,T) \times S^1$ the following equation holds
\begin{equation}\label{eq:loceb}
 \big(W(u) + \frac{1}{2} v^2 + \frac{\gamma}{2} (u_x)^2 \big)_t + \big(-W'(u)v - \mu v v_x + \gamma v u_{xx} -\gamma v_x u_x\big)_x + \mu(v_x)^2 =0.
\end{equation}
Upon integration this implies
\[ \ddt \int_{S^1} W(u) + \frac{1}{2} v^2 + \frac{\gamma}{2} (u_x)^2  \dx = - \mu \int_{S^1} (v_x)^2 \dx \leq 0.\]
\end{lemma}

\bB{}
Equation \eqref{eq:loceb} is obtained by multiplying \eqref{eq:loc}$_1$ by $W'(u) - \gamma u_{xx}$ and \eqref{eq:loc}$_2$ by $v$
and adding both equations.
\eB

Let us define the following non-local surface energy functional
\begin{equation}\label{f1}
 F_\eps: L^2(S^1) \rightarrow [0,\infty) ,\quad w \mapsto  \frac{1}{4\eps^2} \int_{S^1}\int_{S^1} \phi_\eps(x-y)(w(y) - w(x))^2 \dx \dy,  
\end{equation}
where due to the periodic boundary conditions $x-y$ is to be understood as $\dist ( y, \{ x, 1+x, x-1\}).$
\begin{remark}[Reformulation of the energy]
Using the periodicity and the symmetry of $\phi$ we may rewrite $F_\eps$ as follows:
For every $w \in L^2(S^1)$ the following holds
\begin{equation}\label{f2}
 \begin{split}
   F_\eps[w]&= \frac{1}{4\eps^2} \int_{S^1}\int_{S^1} \phi_\eps(x-y)(w(y)^2 -2 w(y) w(x) + w(x)^2) \dx \dy\\
&= \frac{1}{2\eps^2} \int_{S^1}\int_{S^1} \phi_\eps(x-y)(w(y)^2 - w(y) w(x)) \dx \dy\\
&= \frac{1}{2\eps^2} \int_{S^1} w(x)^2 \dx - \frac{1}{2\eps^2} \int_{S^1}\int_{S^1} \phi_\eps(x-y)w(y) w(x) \dx \dy\\
&=- \frac{1}{2} \int_{S^1} w L_\eps[w] \dx.
 \end{split}
\end{equation}
\end{remark}

\begin{lemma}[Energy balance of the non-local model]\label{lem:noneb}
 Let $T,\mu,\eps>0$ and $\phi$ satisfying \eqref{prop:phi} be given. Let $(u^\eps,v^\eps)$ be a strong solution of \eqref{eq:non} in $(0,T)\times S^1$
 in the sense of Lemma \ref{thrm:non}.
Then, for $t \in (0,T)$ the following equation holds
\begin{equation}\label{eq:noneb}
\ddt \Big(  \int_{S^1} (W(u^\eps) + \frac{1}{2} (v^\eps)^2 ) \dx + F_\eps[u^\eps]\Big) = - \int_{S^1} \mu (v_x^\eps)^2 \dx \leq 0.
\end{equation}
\end{lemma}

\bB{}
To obtain equation \eqref{eq:noneb} we multiply \eqref{eq:loc}$_1$ by $W'(u^\eps) - L_\eps [u^\eps]$ and \eqref{eq:loc}$_2$ by $v^\eps$
and add both equations.
Equation \eqref{eq:noneb} follows, due to the boundary conditions upon noting that
\begin{multline}
 -\int_{S^1} u^\eps_t L_\eps [u^\eps]\dx = \frac{1}{\eps^2}\int_{S^1} u^\eps(x) u^\eps_t(x) \dx - \frac{1}{\eps^2}\int_{S^1} \int_{S^1} \phi_\eps(x-y) u^\eps(y) u^\eps_t(x) \dy \dx \\
=  \frac{1}{\eps^2}\int_{S^1} u^\eps(x) u^\eps_t(x) \dx - \frac{1}{2\eps^2} \int_{S^1} \int_{S^1} \phi_\eps(x-y) \big(u^\eps(y) u^\eps_t(x) + u^\eps_t(y) u^\eps(x)\big) \dy \dx
 = \ddt F_\eps[u^\eps],
\end{multline}
because of \eqref{f2} and the symmetry of $\phi$.
\eB

Note that there is a local energy balance for \eqref{eq:loc}, which implies a global version upon integration, but only a global energy balance for \eqref{eq:non}.

\section{Well-posedness and properties of the non-local energy}\label{sec:wp}
In order to carry out our convergence analysis we need to establish the existence of strong solutions to \eqref{eq:loc} and \eqref{eq:non}.
We make use of the standard $L^p$ Lebesgue space notation. By $W^{k,p}$ we denote the Sobolev space of functions with $k$ weak derivatives in $L^p$
and $H^k:=W^{k,2}.$

\subsection{Well-posedness}
For simplicity we complement \eqref{eq:loc} and \eqref{eq:non} with identical initial data 
\begin{equation}\label{ic}
 u(0,\cdot)=u^\eps(0,\cdot)=u_0,\quad  v(0,\cdot)=v^\eps(0,\cdot)=v_0
\end{equation}
for given functions $u_0,v_0 : S^1 \rightarrow \setR.$
We will see later on that it would be sufficient to impose initial data converging to each other in $H^1(S^1) \times L^2(S^1)$ for $\eps \rightarrow 0$, see equation \eqref{eq:t4},  and being sufficiently regular.
In the sequel we choose $u_0 \in H^3(S^1)$.
Let us adopt the convention that for any function space subscript $m$ denotes the  subspace of functions of vanishing mean.
Concerning the well-posedness of \eqref{eq:loc} we cite the following  result from \cite{Gie14}

\begin{lemma}[Well-posedness of \eqref{eq:loc}]\label{thrm:to}
Let initial data $u_0 \in H^3_m(S^1), \, v_0 \in H^2_m(S^1)$ 
 and $T,\mu,\gamma>0$ be given. Then, the problem \eqref{eq:loc}, \eqref{ic} has a unique strong solution
\[ (u,v) \in \big(C^0([0,T],H^3_m(S^1))\cap  C^1((0,T),H^1_m(S^1))\big) \times \big(C^0([0,T],H^2_m(S^1))\cap  C^1((0,T),L^2_m(S^1))\big).\]
\end{lemma}

The existence of strong solutions to \eqref{eq:non} follows from standard semi-group theory.
Our first step is to prove a bound for the non local energy.

\begin{lemma}[$H^1$-stability of non-local energy]\label{ub}
 It exists some $C>0$ independent of $\eps$ such that 
\[ F_\eps[w] \leq C |w|_{H^1(S^1)}^2 \quad \forall \ w \in H^1(S^1).\]
\end{lemma}

\bB{}
Using the definition of $F_\eps$ and Jensen's inequality, we have
\begin{equation}
 \begin{split}
  4F_\eps[w] &= \frac{1}{\eps^2} \int_{S^1} \int_{S^1} \phi_\eps(x-y) (w(y)-w(x))^2 \dx \dy\\
&=\frac{1}{\eps^2} \int_{S^1} \int_{-1}^{1} \phi(s) \Big(\int_{0}^{\eps s} w_x(x + z) \operatorname{d} z\Big)^2 \operatorname{d} s \dx \\
&\leq \int_{S^1} \int_{-1}^{1} \phi(s) \frac{|s|}{\eps} \int_{-|\eps s|}^{|\eps s|} w_x^2(x+ z)  \operatorname{d} z \operatorname{d} s\dx \\
&= \int_{S^1} \int_{-1}^{1} \phi(s) |s| \int_{-|s|}^{|s|} w_x^2(x+\eps \tilde z)  \operatorname{d}\tilde z \operatorname{d} s\dx \\
&=  \int_{-1}^{1}\int_{|s| > |\tilde z|}\int_{S^1} \phi(s) |s|  w_x^2(x+\eps \tilde z) \dx \operatorname{d} s \operatorname{d}\tilde z\\
&\leq  \int_{-1}^{1}\int_{|s| > |\tilde z|} \phi(s) |s|  \operatorname{d} s\operatorname{d}\tilde z\, |w|_{H^1(S^1)}^2 \leq 2 |w|_{H^1(S^1)}^2,
 \end{split}
\end{equation}
which implies the assertion of the lemma.
\eB

\begin{remark}[Regularity of initial data]
 While we will assume $u_0 \in H^3_m(S^1)$ in our convergence analysis, as this is required for the existence of strong solutions to \eqref{eq:loc}, 
we study the well-posedness of \eqref{eq:non} under the weaker assumption $u_0 \in H^1_m(S^1).$
\end{remark}

We continue with a technical lemma which is similar to a result in \cite{ERV13}:

\begin{lemma}[A priori bounds for \eqref{eq:non}]\label{lem:erv}
Provided the assumptions of Lemma \ref{lem:noneb} hold, 
 there is a constant $C>0,$ independent of $\eps,$ such that
$ |u^\eps|_{L^\infty(0,T;H^1(S^1))} <C.$
\end{lemma}
\bB{}
Multiplying \eqref{eq:non}$_2$ by $u_x^\eps$ and integrating in space, we obtain
\begin{equation}
 \int_{S^1} v_t^\eps u_x^\eps - W''(u^\eps)(u_x^\eps)^2 \dx = \int_{S^1} \frac{\mu}{2} ((u_x^\eps)^2)_t - L_\eps[u^\eps]_x u_x^\eps\dx,
\end{equation}
as $v_{xx}^\eps=u_{xt}^\eps.$ Upon  using $L_\eps[u^\eps]_x =L_\eps[u_x^\eps]$ we find
\begin{multline}\label{erv0}
 \int_{S^1}  ((u_x^\eps)^2)_t \dx \leq \frac{2}{\mu}\int_{S^1} (v^\eps u_x^\eps)_t - v^\eps v_{xx}^\eps - W''(u^\eps)(u_x^\eps)^2  + L_\eps[u^\eps_x]u^\eps_x \dx\\
\leq \frac{2}{\mu}\int_{S^1} (v^\eps u_x^\eps)_t - v^\eps v_{xx}^\eps - W''(u^\eps)(u_x^\eps)^2   \dx,
\end{multline}
due to \eqref{f2}. Integrating \eqref{erv0} in time we obtain
\begin{multline}\label{erv1}
 \int_{S^1}  (u_x^\eps)^2(t,\cdot) \dx 
\leq    \frac{2}{\mu}\int_{S^1} \frac{1}{2\delta} (v^\eps)^2(t,\cdot) + \frac{\delta}{2} (u_x^\eps)^2(t,\cdot)
 + (v^\eps)^2(0,\cdot) +  (u_x^\eps)^2(0,\cdot) \dx \\+ \int_{S^1}  (u_x^\eps)^2(0,\cdot) \dx+\frac{2}{\mu} \int_0^t \int_{S^1} (v_x^\eps)^2  - W''(u^\eps)(u_x^\eps)^2 \dx \operatorname{d}s
\end{multline}
for any $\delta >0.$
Choosing $\delta=\tfrac{\mu}{2}$ and the energy dissipation equality \eqref{eq:noneb} we see that \eqref{erv1} together with Lemma \ref{ub}
 implies existence of a constant $C_0 >0$ only depending on 
$u_0,\, v_0$ and $\mu$ such that
\begin{equation}\label{erv2}
 \int_{S^1}  (u_x^\eps)^2(t,\cdot) \dx \leq  C_0 -  \frac{4}{\mu}\int_0^t \int_{S^1}   W''(u^\eps)(u_x^\eps)^2 \dx \operatorname{d}s 
\leq C_0 +  \frac{4\bar W}{\mu}\int_0^t \int_{S^1}    (u_x^\eps)^2 \dx \operatorname{d}s ,
\end{equation}
because of \eqref{ass:W}.
The assertion of the Lemma follows from \eqref{erv2} upon using  Gronwall's inequality.
\eB

\begin{lemma}[Well-posedness of \eqref{eq:non}]\label{thrm:non}
Let initial data $u_0 \in H^1_m(S^1), \, v_0 \in H^2_m(S^1)$ 
 and $T,\mu,\eps>0$ be given. Then, the problem \eqref{eq:non}, \eqref{ic} has a unique strong solution $(u^\eps,v^\eps)$ with
\[
\begin{split} u^\eps&\in 
 W^{1,\infty}(0,T;L^2_m(S^1)) \cap H^1(0,T;H^1_m(S^1)),\\
v^\eps &\in L^2(0,T;H^2_m(S^1)) \cap L^\infty(0,T;H^1_m(S^1)) \cap  H^1(0,T;L^2_m(S^1)).
 \end{split}
\]
\end{lemma}
\bB{}
Let us rewrite \eqref{eq:non} using  the unique primitive $\tau^\eps$ of $u^\eps$ (with respect to $x$) with vanishing mean.
Note that for every $t$ the primitive $\tau^\eps(t,\cdot)$ is periodic as $u^\eps(t, \cdot)$ has vanishing mean.
We will study the following system which is equivalent to \eqref{eq:loc}
\begin{equation}\label{erv3}
 \begin{split}
  \tau_t^\eps &= v^\eps\\ v_t^\eps - W'(\tau_x^\eps)_x&= \mu v_{xx}^\eps  -\frac{1}{\eps^2} \phi_\eps*\tau_{xx}^\eps +\frac{1}{\eps^2} \tau_{xx}^\eps
 \end{split}
\end{equation}
with initial data $(\tau_0,v_0)$ where $\tau_0$ is the unique primitive of $u_0$ in $H^2_m(S^1).$
We will show that \eqref{erv3} has strong solutions on arbitrarily large time intervals.
The proof of this consists of four steps:

Step 1: We consider the reduced problem
\begin{equation}\label{erv4}
 \begin{split}
  \tau_t^\eps &= v^\eps\\ v_t^\eps&= \mu v_{xx}^\eps +\frac{1}{\eps^2} \tau_{xx}^\eps . 
 \end{split}
\end{equation}
Equation \eqref{erv4} can be written in abstract form as 
\begin{equation}\label{erv5}
 \ddt z^\eps(t)=A z^\eps(t) \quad \text{with } z^\eps:=\begin{pmatrix} \tau^\eps \\ v^\eps\end{pmatrix} \ \text{and } \ A:= \begin{pmatrix}
                                                                                                  0 & \Id \\ \eps^{-2}\del_{xx} & \mu \del_{xx}
                                                                                                 \end{pmatrix}
\end{equation}
with $D(A)= H^2_m(S^1) \times H^2_m(S^1).$
It is shown in \cite{EN00} that $A$ generates a strongly continuous semi-group on $Y=H^2_m(S^1) \times L^2_m(S^1).$

Step 2: In this step, we consider the system 
\begin{equation}\label{erv6}
 \begin{split}
  \tau_t^\eps &= v^\eps\\ v_t^\eps&= \mu v_{xx}^\eps  -\frac{1}{\eps^2} \phi_\eps*\tau_{xx}^\eps+\frac{1}{\eps^2} \tau_{xx}^\eps. 
 \end{split}
\end{equation}
It can be written as
\begin{equation}\label{erv7} \ddt z^\eps(t)=A z^\eps(t) + B z^\eps(t)\end{equation}
using the notation from \eqref{erv5} and a bounded operator $B$ on $Y.$
Thus, \cite[Theorem III.1.3]{EN00} implies that \eqref{erv6} admits strong solutions.

Step 3: Using the notation from \eqref{erv7} we consider
\begin{equation}\label{erv8} \ddt z^\eps(t)=A z^\eps(t) + B z^\eps(t) + F(z^\eps(t))\end{equation}
with
\[ F(z^\eps) := \begin{pmatrix}
              0\\ W'(\tau_x^\eps)_x
             \end{pmatrix}.\]
As observed in \cite[e.g.]{Gie14} the map $F:Y \rightarrow Y$ is locally Lipschitz and, thus, by Banach's fixed point Theorem
\eqref{erv8} has a unique mild solution for sufficiently small times. As the initial data are in $D(A)$ this mild solution is, in fact, a 
strong solution.

Step 4: We will show that the local-in-time solution in Step 3 can be extended to arbitrary times.
To this end we show that the solution does not blow up in finite time.
We know from Lemmas \ref{lem:noneb} and \ref{lem:erv} that $\|v^\eps(t,\cdot)\|_{L^2(S^1)}$ and $\|\tau_{xx}^\eps(t,\cdot)\|_{L^2(S^1)}$ remain bounded for arbitrary times.
Therefore, we only need to control $v_x^\eps$ and $v_{xx}^\eps.$
Let us note, that \eqref{eq:non}$_2$ can be written as an inhomogeneous heat equation
\begin{equation}\label{erv9}
 v_t^\eps -\mu v_{xx}^\eps =g^\eps:= W'(\tau_x^\eps)_x - L_\eps[\tau_x^\eps]_x
\end{equation}
and the right hand side is controlled by 
\[\| g^\eps \|_{L^2(0,T;L^2(S^1))} \leq C \|\tau_{xx}^\eps\|_{L^2(0,T;L^2(S^1))},\]
for some  constant $C>0$, 
due to the embedding of $H^1(S^1)$ in $C^0(S^1)$ and the regularity of $W$.
Standard elliptic regularity theory, i.e. squaring both sides of \eqref{erv9} and integrating, yields for all $ T>0$
\begin{multline}
 \| v_t^\eps\|_{L^2(0,T;L^2(S^1))}^2 + \mu \| v_x^\eps\|_{L^\infty(0,T;L^2(S^1))}^2  +\mu^2 \| v_{xx}^\eps\|_{L^2(0,T;L^2(S^1))}^2\\ 
\leq \mu \| v_x^\eps(0,\cdot)\|_{L^2(S^1)}^2 + \| g^\eps\|_{L^2(0,T;L^2(S^1))}^2
\leq \mu | v_0 |_{H^1(S^1)}^2 + C \| \tau_{xx}^\eps\|_{L^2(0,T;L^2(S^1))}^2.
\end{multline}
Thus, $v$ has the asserted regularity.
Equation \eqref{erv3} and the regularity of $v^\eps$ imply 
\[ \tau^\eps \in W^{1,\infty}(0,T;H^1_m(S^1)) \cap H^1(0,T;H^2_m(S^1))\]
such that
\[u^\eps\in  W^{1,\infty}(0,T;L^2_m(S^1)) \cap H^1(0,T;H^1_m(S^1)),\]
which is the desired regularity for $u.$
\eB


\section{The non-local to local limit}\label{sec:limit}
Let us define the following variant of the relative entropy between a solution $(u,v)$ of \eqref{eq:loc} and  a solution $(u^\eps,v^\eps)$ 
of \eqref{eq:non}:
\begin{multline}\label{def:re}
  \eta_\eps(t) := F_\eps [u^\eps] + \int_{S^1} W(u^\eps) + \frac{1}{2} (v^\eps)^2 - W(u) - \frac{1}{2} v^2 - \frac{\gamma}{2} (u_x)^2\\
 -W'(u)(u^\eps-u) - v(v^\eps-v) + \gamma u_x u_x + L_\eps[u]u^\eps \dx,
 \end{multline}
where it is  understood that all the functions on the right hand side are to be evaluated at time $t.$
It is important to note that $\eta_\eps$ is not convex with respect to $(u - u^\eps,v-v^\eps)$ mainly due to the $W$ terms.
 Therefore,  we  like to introduce a modified relative entropy without the $W$ terms, but including $\| u^\eps-u\|^2_{L^2(S^1)}:$
\begin{equation}\label{def:rre}
   \eta_\eps^{\textrm{M}}(t) := F_\eps [u^\eps] + \int_{S^1}  \frac{1}{2} (v^\eps)^2 - \frac{1}{2} v^2 - \frac{\gamma}{2} (u_x)^2
 - v(v^\eps-v) + \gamma u_x u_x + L_\eps[u]u^\eps + \frac{1}{2} ( u^\eps-u)^2 \dx.
\end{equation}
Let us note that several terms in $\eta_\eps$ and $ \eta_\eps^{\textrm{M}}$ cancel out such that
\begin{equation}\label{eq:re}
 \begin{split}
  \eta_\eps(t) =& F_\eps [u^\eps(t,\cdot)] + \int_{S^1} W(u^\eps(t,\cdot)) - W(u(t,\cdot))-W'(u(t,\cdot))(u^\eps(t,\cdot)-u(t,\cdot))
\\&+ \frac{1}{2} (v^\eps(t,\cdot) - v(t,\cdot))^2
  + \frac{\gamma}{2} (u_x(t,\cdot))^2 + L_\eps[u(t,\cdot)] u^\eps(t,\cdot) \dx,\\
   \eta_\eps^{\textrm{M}}(t) &= F_\eps [u^\eps(t,\cdot)] + \int_{S^1}  \frac{1}{2} (v^\eps(t,\cdot)- v(t,\cdot))^2   + \frac{1}{2} ( u^\eps(t,\cdot)-u(t,\cdot))^2\\
& + \frac{\gamma}{2} (u_x(t,\cdot))^2
  + L_\eps[u(t,\cdot)] u^\eps(t,\cdot)\dx.
 \end{split}
\end{equation}

\begin{lemma}[Relative entropy rate]\label{lem:re}
 Let $u_0 \in H^3_m(S^1),$ $v_0 \in H^2_m(S^1)$, $T,\mu,\gamma,\eps>0$ and $\phi$ satisfying \eqref{prop:phi} be given.
Then, the strong solution $(u,v)$ of \eqref{eq:loc}, \eqref{ic} and the strong solution $(u^\eps,v^\eps)$ of  \eqref{eq:non}, \eqref{ic}
satisfy
\begin{equation}\label{eq:rel}
 \ddt \eta_\eps = \int_{S^1} v_x \big(W'(u^\eps) - W'(u) - W''(u)(u^\eps-u)
\big) - \mu (v_x - v^\eps_x)^2 + \gamma v^\eps u_{xxx} + v^\eps_x L_\eps[u] \dx.
\end{equation}
\end{lemma}

\bB{}
Using the same identities for the time derivative of $F[u^\eps]$ as in the proof of Lemma \ref{lem:noneb} we find
\begin{multline}\label{eq:ree1}
  \ddt \eta_\eps
= \int_{S^1} -L_\eps[u^\eps]u^\eps_t + v^\eps v^\eps_t + W'(u^\eps)u^\eps_t - W'(u) u_t + v v_t + \gamma u_x u_{xt} - W''(u)u_t (u^\eps - u) \\- W'(u) u^\eps_t + W'(u)u_t
-v v^\eps_t - v_t v^\eps + L_\eps[u_t] u^\eps + L_\eps[u] u^\eps_t\dx.
 \end{multline}
Replacing the time derivatives in \eqref{eq:ree1} using the evolution equations \eqref{eq:loc}, \eqref{eq:non} we infer
\begin{multline}\label{eq:ree2}
  \ddt \eta_\eps
= \int_{S^1} -L_\eps[u^\eps]v^\eps_x + v^\eps (W'(u^\eps)_x + \mu v^\eps_{xx} - L_\eps[u^\eps]_x) + W'(u^\eps) v^\eps_x  + v (W'(u)_x + \mu v_{xx} - \gamma u_{xxx})\\
 + \gamma u_x v_{xx} - W''(u)v_x (u^\eps - u) - W'(u) v^\eps_x
-v(W'(u^\eps)_x + \mu v^\eps_{xx} - L_\eps[u^\eps]_x)\\ - v^\eps(W'(u)_x + \mu v_{xx} - \gamma u_{xxx}) + L_\eps[v_x] u^\eps + L_\eps[u] v^\eps_x\dx .
 \end{multline}
Rewriting the right hand side of \eqref{eq:ree2} we obtain
\begin{multline}\label{eq:ree3}
  \ddt \eta_\eps
= \int_{S^1} \Big( -L_\eps[u^\eps] v^\eps + (v^\eps - v) (W'(u^\eps)-W'(u)) - \gamma u_{xx}   v + \gamma u_x   v_x \Big)_x\dx\\
+\int_{S^1} v_x \big(W'(u^\eps) - W'(u) - W''(u)(u^\eps - u) \big) +\mu(v^\eps -v)_{xx} (v^\eps - v)\\
 + v L_\eps[u^\eps]_x +\gamma v^\eps u_{xxx} + L_\eps[v_x] u^\eps + L_\eps[u] v^\eps_x \dx.
 \end{multline}
Due to the symmetry of $\phi$ and Fubini's Theorem, we have, for almost every $t \in [0,T],$
\begin{multline}\label{eq:ree4}
  \int_{S^1} v(t,\cdot) L_\eps[u^\eps(t,\cdot)]_x  \dx 
=\frac{1}{\eps^2} \int_{S^1} \Big( \int_{S^1}  \phi_\eps (x-y) u_x^\eps(t,y)\dy - u_x^\eps(t,x)\Big) v(t,x)\dx \\
=\frac{1}{\eps^2} \int_{S^1} \Big( \int_{S^1}  \phi_\eps (x-y) v(t,x)\dx - v(t,y)\Big) u_x^\eps(t,y)\dy
=   \int_{S^1} L_\eps[v(t,\cdot)] u_x^\eps(t,\cdot)\dx\\
= - \int_{S^1} L_\eps[v(t,\cdot)]_x u^\eps(t,\cdot)\dx
=  -\int_{S^1} L_\eps[v_x(t,\cdot)] u^\eps(t,\cdot)\dx.
\end{multline}
Due to the periodic boundary conditions and \eqref{eq:ree4}, \eqref{eq:ree3} implies the assertion of the Lemma.
\eB

Our next step is to derive an estimate for the rate of $\eta_\eps^{\textrm{M}}.$
It can be obtained from Lemma  \ref{lem:re} analogous to \cite[Cor 3.3]{Gie14}.
We give the proof for completeness.
\begin{lemma}[Estimate of the rate of $\eta_\eps^{\textrm{M}}$]\label{lem:rre}
 Let $u_0 \in H^3_m(S^1),$ $v_0 \in H^2_m(S^1)$, $T,\mu,\gamma,\eps>0$ and $\phi$ satisfying \eqref{prop:phi} be given.
Then, the strong solution $(u,v)$ of \eqref{eq:loc}, \eqref{ic} and the strong solution $(u^\eps,v^\eps)$ of  \eqref{eq:non}, \eqref{ic}
satisfy
\begin{equation}\label{eq:rre}
 \ddt \eta_\eps^{\textrm{M}} \leq \int_{S^1} \frac{1}{2\mu} \big(W'(u^\eps) - W'(u) \big)^2   + \frac{1}{2\mu} \big(u^\eps - u \big)^2 + \gamma v^\eps u_{xxx} + v^\eps_x L_\eps[u] \dx.
\end{equation}
\end{lemma}

\bB{}
We note that
\begin{equation}\label{rrer1}
 \eta_\eps = \eta_\eps^{\textrm{M}} + \int_{S^1} W(u^\eps) - W(u) - W'(u)(u^\eps - u) -   \frac{1}{2} ( u^\eps-u)^2 \dx
\end{equation}
and
\begin{equation}\label{rrer2}
\begin{split}
 \frac{\operatorname{d}}{\dt} \int_{S^1} W(u^\eps) - W(u) - W'(u)(u^\eps - u)\dx&= \int_{S^1} v^\eps_x (W'(u^\eps)-W'(u)) - v_x W''(u) (u^\eps - u) \dx,\\
- \frac{1}{2} \frac{\operatorname{d}}{\dt} \int_{S^1} (u^\eps-u)^2 \dx&= -\int_{S^1} (u^\eps-u)(v^\eps_x-v_x) \dx.
\end{split} 
\end{equation}
Combining \eqref{rrer1} and \eqref{rrer2} with \eqref{eq:rel} we find 
\begin{equation}
 \ddt \eta_\eps^{\textrm{M}} = \int_{S^1} (v_x- v^\eps_x) \big(W'(u^\eps) - W'(u) 
\big) + (u^\eps-u)(v^\eps_x-v_x) -\mu (v_x - v^\eps_x)^2 + \gamma v^\eps u_{xxx} + v^\eps_x L_\eps[u] \dx.
\end{equation}
We obtain the assertion of the Lemma upon using Young's inequality.
\eB

Our next step is to estimate the last two summands on the right hand side of \eqref{eq:rre}.
They stem from the difference between the local and non-local regularization.
In most formal arguments linking non-local to local models it was argued that $L_\eps[u^\eps] \rightarrow \gamma u^\eps_{xx}$
for $\eps \rightarrow 0.$
However, in general, $u^\eps_{xx}$ is not well-defined such that we need to modify this approach. 
The main idea is to analyze $u$ instead of $u^\eps.$
\begin{lemma}[Relative entropy increase by difference of regularizations]\label{lem:error}
 Let the assumptions of Lemma \ref{lem:rre} be satisfied. 
Then,
\[\int_0^T \Big| \int_{S^1}  \gamma v^\eps u_{xxx} + v^\eps_x L_\eps[u]   \dx \Big|\dt = \mathcal{O}(\sqrt{\eps}).
\]
\end{lemma}
\bB{}
We have
\begin{multline}\label{fs}
\int_0^T\Big| \int_{S^1} v^\eps_x(t,\cdot) L_\eps[u(t,\cdot)]+  \gamma v^\eps(t,\cdot) u_{xxx}(t,\cdot)   \dx \Big|\dt \\
= \int_0^T\Big| \int_{S^1}  (L_\eps[u(t,\cdot)]-  \gamma u_{xx}(t,\cdot)) v^\eps_x(t,\cdot)   \dx \Big|\dt\\
\leq \sqrt{T}  \| L_\eps[u]-  \gamma u_{xx}\|_{L^2([0,T]\times S^1)} \|v^\eps_x\|_{L^2([0,T]\times S^1)}.
\end{multline}
We know from Lemma \ref{lem:noneb} that $\|v^\eps_x\|_{L^2([0,T]\times S^1)}$ is  bounded in terms of the initial energy and $\mu$.
We are going to show 
$\| L_\eps[u]-  \gamma u_{xx}\|_{L^2(0,T;L^2(S^1))} =\mathcal{O}(\sqrt{\eps})$ for $\eps \rightarrow 0.$
To this end, let us fix some $(t,x) \in [0,T] \times S^1.$ We have
\begin{multline}\label{eq:v1}
 \int_{S^1} \phi_\eps(x-y) (u(t,y) - u(t,x)) \dy = \int_{-1}^1 \phi(s)  (u(t,x+s\eps) - u(t,x)) \operatorname{d}s\\
=\int_{-1}^1 \phi(s) \eps \int_0^s u_x(t,x+ \eps z)  \operatorname{d}z\operatorname{d}s
= \int_{-1}^1 \phi(s) \eps \int_0^s \Big(u_x(t,x) + \int_0^{\eps z}  u_{xx}(t,x+ a) \operatorname{d}a \Big)  \operatorname{d}z\operatorname{d}s.
\end{multline}
Due to the symmetry of $\phi,$ equation \eqref{eq:v1} implies 
\begin{equation}\label{eq:v2}
  \frac{1}{\eps^2}\int_{S^1} \phi_\eps(x-y) (u(t,y) - u(t,x)) \dy 
= \int_{-1}^1 \phi(s) \frac{1}{\eps} \int_0^s  \int_0^{\eps z}  u_{xx}(t,x+ a)
 \operatorname{d}a  \operatorname{d}z\operatorname{d}s.
\end{equation}
As $u_{xxx} \in C^0([0,T],L^2( S^1)),$ we have, for all $(t,x) \in [0,T] \times S^1,$ 
 \begin{equation}\label{**}
  \begin{split}
 & \Big| \frac{1}{\eps} \int_0^{\eps z} u_{xx}(t,x+a) \operatorname{d} a - z u_{xx} (t,x) \Big|\\
&\leq\frac{1}{\eps} \int_{-\eps|z|}^{\eps |z|} \Big|\int_0^a |u_{xxx}(t,x+\xi)| \operatorname{d} \xi\Big|\operatorname{d} a  
= \frac{1}{\eps} \int_{-\eps|z|}^{\eps |z|} \int_{|\xi|< |a| < \eps |z|} |u_{xxx}(t,x+\xi)|  \operatorname{d} a\operatorname{d} \xi\\
&\leq \frac{2}{\eps} \int_{-\eps|z|}^{\eps |z|}  |u_{xxx}(t,x+\xi)| \, \Big| \eps |z| - |\xi| \Big| \operatorname{d} \xi
\leq 2|z|\int_{-\eps|z|}^{\eps |z|}  |u_{xxx}(t,x+\xi)|\operatorname{d} \xi
 \leq 4\sqrt{\eps |z|^3} \| u_{xxx}(t,\cdot)\|_{L^2(S^1)}
  \end{split}
 \end{equation}
and $\| u_{xxx}(t,\cdot)\|_{L^2(S^1)}$ is uniformly bounded in time.
Thus, for all $(t,x)$
\begin{equation}\label{eq:v3}
\begin{split}
   &\frac{1}{\eps^2}\int_{S^1} \phi_\eps(x-y) (u(t,y) - u(t,x)) \dy \\
&= \int_{-1}^1 \phi(s) \int_0^s z u_{xx}(t,x) + \Big( \frac{1}{\eps} \int_0^{\eps z} u_{xx}(t,x+a) \operatorname{d} a - z u_{xx} (t,x) \Big)\operatorname{d} z \operatorname{d} s\\
&= \int_{-1}^1 \phi(s) \frac{s^2}{2}  u_{xx}(t,x) \operatorname{d} s + \mathcal{O}(\sqrt{\eps}),
\end{split}
\end{equation}
uniformly in space and time, which implies 
\begin{equation}\label{ls}
  \| L_\eps[u]-  \gamma u_{xx}\|_{L^2(0,T;L^2(S^1))} \leq T  \| L_\eps[u]-  \gamma u_{xx}\|_{L^\infty((0,T)\times S^1)} \leq  \mathcal{O}(\sqrt{\eps}).
\end{equation}
The assertion of the Lemma follows upon combining \eqref{fs} and \eqref{ls}.
\eB

Our next Lemma studies the surface energy terms in $\eta_\eps^{\textrm{M}}$ in more detail.
\begin{lemma}[Surface terms in  $\eta_\eps^{\textrm{M}}$]\label{lem:convex}
 Let the assumptions of Lemma \ref{lem:rre} be satisfied. Then,
\[  F_\eps [u^\eps(t,\cdot)] + \int_{S^1}    \frac{\gamma}{2} (u_x(t,\cdot))^2
  + L_\eps[u(t,\cdot)] u^\eps(t,\cdot) \dx -  F_\eps [u^\eps(t,\cdot) -u(t,\cdot)] = \mathcal{O}(\eps)
\]
for $\eps \rightarrow 0,$ uniformly in time.
\end{lemma}
\bB{}
As in the proof of Lemma \ref{lem:error}, we do not study the terms in  $\eta_\eps^{\textrm{M}}$ which are related to $u^\eps$ but those related to $u$.
 In this way we can make use of the 
higher regularity of $u.$
Let us note that
\begin{equation}\label{conv1}
\begin{split}
  F_\eps[u(t,\cdot)]&=\frac{1}{4\eps^2} \int_{S^1}\int_{S^1} \phi_\eps(x-y) \Big( \int_x^y u_x(t,z)\operatorname{d}z \Big)^2 \dy \dx\\
        &= \frac{1}{4\eps^2}\int_{S^1}\int_{S^1} \phi_\eps(x-y) \Big(\int_x^y \Big( u_x(t,x) + \int_x^z u_{xx}(t,a)\operatorname{d}a\Big)\operatorname{d}z \Big)^2\dy \dx\\
         &= \int_{S^1}\frac{\gamma}{2} (u_x)^2(t,\cdot) \dx +  \|u(t,\cdot)\|_{C^2(S^1)} \eps
\end{split}
\end{equation}
because $u \in C^0([0,T],H_m^3(S^1)) \subset C^0([0,T],C^2(S^1)).$
Moreover, due to the symmetry of $\phi,$
\begin{equation}\label{conv2}
\begin{split}
& \frac{1}{2\eps^2} \int_{S^1}\int_{S^1} \phi_\eps(x-y) \big( u(t,x) - u(t,y)\big) \big( u^\eps(t,x) - u^\eps(t,y)\big) \dy \dx\\
&= \frac{1}{\eps^2}\int_{S^1}\int_{S^1} \phi_\eps(x-y) \big( u(t,x) - u(t,y)\big)  u^\eps(t,x)\dy \dx\\
&= -\int_{S^1} L_\eps[u(t,\cdot)] u^\eps(t,\cdot) \dx.
\end{split}
\end{equation}
Equations \eqref{conv1} and \eqref{conv2} imply
\begin{multline}
 F_\eps[u^\eps(t,\cdot)] + \int_{S^1}    \frac{\gamma}{2} (u_x(t,\cdot))^2 + L_\eps[u(t,\cdot)]u^ \eps(t,\cdot) \dx\\
= F_\eps[u^\eps(t,\cdot)] + F_\eps[u(t,\cdot) ] - \frac{1}{2\eps^2} \int_{S^1}\int_{S^1} \phi_\eps(x-y) \big( u(t,x) - u(t,y)\big) \big( u^\eps(t,x) - u^\eps(t,y)\big)\dx \dy +\mathcal{O}(\eps)\\
=F_\eps[u^\eps(t,\cdot) - u(t,\cdot)] +\mathcal{O}(\eps),
\end{multline}
which is the assertion of the lemma.
\eB

\begin{theorem}[Model convergence]\label{thrm}
 Let $u_0 \in H^3_m(S^1),$ $v_0 \in H^2_m(S^1)$, $T,\mu,\gamma >0$ and $\phi$ satisfying \eqref{prop:phi} be given.
Then, for $\eps>0$ sufficiently small, the strong solution $(u,v)$ of \eqref{eq:loc}, \eqref{ic} and the strong solution $(u^\eps,v^\eps)$ of  \eqref{eq:non}, \eqref{ic} 
satisfy 
\begin{equation}\label{*} \| u^\eps - u\|_{L^\infty(0,T;L^2(S^1))}^2 +\| v^\eps - v\|_{L^\infty(0,T;L^2(S^1))}^2  \leq C \sqrt{\eps} \end{equation}
with $C>0$ depending on $u_0,v_0,\mu,\gamma,T,\phi.$
\end{theorem}

\bB{}
Integrating the assertion of Lemma \ref{lem:rre} in time and using Lemma \ref{lem:error} we find
\begin{equation}\label{eq:t1}
  \eta_\eps^{\textrm{M}}(t) - \eta_\eps^{\textrm{M}}(0) \leq \int_0^t \int_{S^1} \frac{1}{2\mu} \big(W'(u^\eps) - W'(u) \big)^2 + \frac{1}{2\mu} \big(u^\eps-u \big)^2 \dx \operatorname{d}s + \mathcal{O}(\sqrt{\eps}).
\end{equation}
As $H^1(S^1)$ is continuously embedded in $C^0(S^1)$ we know that $\|u^\eps\|_{L^\infty((0,T)\times S^1)}$, $\|u\|_{L^\infty((0,T)\times S^1)}$ are bounded independent of $\eps,$ cf. Lemma \ref{lem:erv},
 such that the regularity of $W$ implies
\begin{equation}\label{eq:t2}
  \eta_\eps^{\textrm{M}}(t) - \eta_\eps^{\textrm{M}}(0)  \leq C \int_0^t \|u^\eps(s, \cdot) - u(s,\cdot) \|_{L^2(S^1)}^2 \operatorname{d}s+ \mathcal{O}(\sqrt{\eps})
\end{equation}
for some constant $C$ independent of $\eps.$
Upon applying Lemma \ref{lem:convex} and using the non-negativity of $F_\eps$ we infer from equation \eqref{eq:t2}
\begin{equation}\label{eq:t3}
\begin{split}
&  \|u^\eps(t,\cdot) - u(t,\cdot)\|_{L^2(S^1)}^2+ \|v^\eps(t,\cdot) - v(t,\cdot)\|_{L^2(S^1)}^2 + \mathcal{O}(\eps)\\
  \leq & C \int_0^t \|u^\eps(s,\cdot) - u(s,\cdot)\|_{L^2(S^1)}^2\operatorname{d}s + \mathcal{O}(\sqrt{\eps})\\
\leq & C \int_0^t \Big( \|u^\eps(s,\cdot) - u(s,\cdot)\|_{L^2(S^1)}^2+ \|v^\eps(s,\cdot) - v(s,\cdot)\|_{L^2(S^1)}^2  + \mathcal{O}(\eps) \Big)\operatorname{d}s + (C+1) \mathcal{O}(\sqrt{\eps}),
\end{split}
\end{equation}
where we have, in particular, used that Lemma  \ref{lem:convex} implies 
\begin{multline}
 \eta_\eps^{\textrm{M}}(0) = F_\eps[u^\eps(0,\cdot)] +  \int_{S^1}\frac{\gamma}{2} |u_x(0,\cdot)|^2 + L_\eps[u(0,\cdot)]u^\eps(0,\cdot) \operatorname{d} x\\
= F_\eps[u^\eps(0,\cdot) - u(0,\cdot)] + \mathcal{O}(\eps)= \mathcal{O}(\eps),
\end{multline}
as the initial data coincide. Using Gronwall's Lemma we obtain
\begin{equation}\label{eq:t4}
 \|u^\eps - u\|_{L^\infty(0,T;L^2(S^1))}^2 +  \|v^\eps - v\|_{L^\infty(0,T;L^2(S^1))}^2  + \mathcal{O}(\eps)  \leq   \mathcal{O}(\sqrt{\eps}) e^{CT} .
\end{equation}
This implies 
\begin{equation}\label{eq:t5}
 \|u^\eps - u\|_{L^\infty(0,T;L^2(S^1))}^2 +  \|v^\eps - v\|_{L^\infty(0,T;L^2(S^1))}^2   \leq   \mathcal{O}(\sqrt{\eps}) e^{CT} +\mathcal{O}(\sqrt{\eps})
\end{equation}
from which the assertion of the Theorem follows.
\eB

\begin{remark}[Multiple space dimensions]
 We expect that the arguments presented here can be formally extended to several space dimensions, as in a similar investigation in \cite{Gie14}. 
However, we do not pursue this approach here as the multi-dimensional, generalized version of the local model \eqref{eq:loc} considered in \cite{Gie14} is physically inadmissible and determining a physically
admissible generalization is the subject of ongoing research
\end{remark}

\begin{remark}[Parameter and time dependence]
 It must be noted that the convergence result derived in Theorem \ref{thrm} depends sensitively on $\mu.$
Thus, the non-viscous ($\mu \rightarrow 0$) limit cannot be studied by the arguments presented in this work.
The dependence on $\phi$ and $\gamma$ is more subtle, and enters mainly in the constant in Lemma \ref{lem:convex} and via the properties of the solution of the local model.
While we have uniform convergence on compact time intervals, the convergence becomes slower for larger times as can be seen from \eqref{eq:t5}.
\end{remark}

\begin{remark}[Convergence rate]
 Note that we are mainly interested in the convergence $(u^\eps,v^\eps) \rightarrow (u,v)$ as such, and the convergence rate is of minor importance.
 However, we like to point out that in case $u_{xxx} \in L^\infty([0,T] \times S^1)$ we could replace $\sqrt{\eps}$ in \eqref{*} by $\eps$, i.e., the 
 convergence rate would double. 
 This can be seen by redoing the last step of \eqref{**} and tracking the power of $\eps$ in the proof of Theorem \ref{thrm}.
\end{remark}

\begin{remark}[Smaller state space]\label{rem:state}
It is possible to extend our results to the case of $W$  only being defined on some interval $I \subset\setR$ and $W(u)$ diverging for $u \rightarrow \del I.$
This is in particular important if we view \eqref{eq:loc} and \eqref{eq:non} as models for liquid vapor flows or longitudinal motions of an elastic bar.
In those cases physics, i.e., the fact that there is no interpenetration of matter, requires $u >0.$
 However, in that case we would need to restrict ourselves to situations in which $u,u^\eps$ are uniformly bounded away from $\del I.$
\end{remark}

Let us conclude by pointing out that the modified relative entropy framework offers an easy proof of continuous dependence of strong 
 solutions to \eqref{eq:non} on their initial data.

\begin{lemma}[Continuous dependence on initial data]\label{Continuous dependence}
 Let $\eps, \mu, T>0$ and $\phi$ satisfying \eqref{prop:phi} for some $\gamma>0$ be given.  Let $(u^\eps, v^\eps), (\tilde u^\eps, \tilde v^\eps)$ 
be strong solutions to \eqref{eq:non} on $[0,T]\times S^1$ corresponding to initial data
$(u_0, v_0), (\tilde u_0, \tilde v_0) \in H^1_m(S^1) \times H^2_m(S^1),$ respectively.
 Then, there exists a constant $C>0$ depending on $(u_0, v_0, \| \tilde u_0\|_{H^1(S^1)}, \eps, \mu,T)$ such that
\begin{equation} \|u^\eps - \tilde u^\eps\|_{L^\infty(0,T;L^2(S^1))} + \|v^\eps - \tilde v^\eps\|_{L^\infty(0,T;L^2(S^1))} \leq C 
\Big( \|u_0 - \tilde u_0\|_{L^2(S^1)} + \|v_0 - \tilde v_0\|_{L^2(S^1)}\Big).
\end{equation}
\end{lemma}

\bB{}
 Firstly we note that $u^\eps, \tilde u^\eps$ are bounded in $L^\infty([0,T] \times S^1)$ by some constant depending on the same quantities as $C$ in the assertion of the Lemma,
by the arguments of the proof of Lemma \ref{lem:erv}.
Analogous to the derivation of Lemma \ref{lem:rre} we can show that
\begin{multline}
 \frac{\operatorname{d}}{\operatorname{d} t} \Big(  \| u^\eps(t,\cdot) - \tilde u^\eps(t,\cdot)\|_{L^2(S^1)}^2 + \| v^\eps(t,\cdot) - \tilde v^\eps(t,\cdot)\|_{L^2(S^1)}^2 + F_\eps[u^\eps(t,\cdot) - \tilde u^\eps(t,\cdot)] \Big) \\
\leq  \frac{1}{2\mu} \| W'(u^\eps(t,\cdot)) - W'(\tilde u^\eps(t,\cdot)) \|_{L^2(S^1)}^2 + \frac{1}{2\mu} \| u^\eps(t,\cdot)-\tilde u^\eps(t,\cdot) \|_{L^2(S^1)}^2
\end{multline}
for $ 0<t <T.$
Due to the uniform bound on $u^\eps, \tilde u^\eps$  in $L^\infty([0,T] \times S^1)$ we can use the regularity of $W$  to infer
\begin{multline}
 \frac{\operatorname{d}}{\operatorname{d} t} \Big(\| u^\eps(t,\cdot) - \tilde u^\eps(t,\cdot)\|_{L^2(S^1)}^2 
+ \| v^\eps(t,\cdot) - \tilde v^\eps(t,\cdot)\|_{L^2(S^1)}^2 + F_\eps[u^\eps(t,\cdot) - \tilde u^\eps(t,\cdot)] \Big) 
\\ \leq  \frac{C}{\mu} \| u^\eps(t,\cdot) - \tilde u^\eps(t,\cdot)\|_{L^2(S^1)}^2,
\end{multline}
such that Gronwall's Lemma and the non-negativity of $F_\eps$ imply
\begin{multline}\label{stab1}
\| u^\eps(t,\cdot) - \tilde u^\eps(t,\cdot)\|_{L^2(S^1)}^2 +  \| v^\eps(t,\cdot) - \tilde v^\eps(t,\cdot)\|_{L^2(S^1)}^2  \\
\leq   \big( \| u^\eps(0,\cdot) - \tilde u^\eps(0,\cdot)\|_{L^2(S^1)}^2  + \| v^\eps(0,\cdot) - \tilde v^\eps(0,\cdot)\|_{L^2(S^1)}^2 
+ F_\eps[u^\eps(0,\cdot) - \tilde u^\eps(0,\cdot)] \big) e^{\frac{C}{\mu}t}.
\end{multline}
It is easy to show that for all $w \in L^2(S^1)$ it holds
\begin{equation}\label{star} F_\eps[w] \leq \frac{1}{\eps^2} \| w\|_{L^2(S^1)}^2.\end{equation}
Using \eqref{star}  we infer the assertion of the Lemma at hand from \eqref{stab1}.
\eB

\thanks{
We gratefully acknowledge that this work  was supported by the 
German Research Foundation (DFG) via SFB TRR 75 `Tropfendynamische Prozesse unter extremen Umgebungsbedingungen'.
}

\bibliographystyle{plain}      
\bibliography{nskbib}
\end{document}